%% file: converse_Bohr_equivalence_thm.tex
\documentclass[10pt]{amsart}

\usepackage[a4paper]{geometry}
\usepackage{stile_en}

\usepackage[pdftex, plainpages = false, pdfpagelabels, 
						 bookmarks,
						 bookmarksopen = true,
						 bookmarksnumbered = true,
						 breaklinks = true,
						 linktocpage,
						 colorlinks = true,
						 linkcolor = blue,
						 urlcolor  = blue,
						 citecolor = red,
						 anchorcolor = green,
						 hyperindex = true,
						 hyperfigures
						 ]{hyperref} 
\usepackage[pdftex]{graphicx}
\usepackage[utf8]{inputenc}
\usepackage[english]{babel}
\usepackage[T1]{fontenc}

\usepackage{type1cm}

\usepackage{amsmath, amsfonts, amssymb}
\usepackage[shortlabels]{enumitem}
\usepackage[all]{xy}

\usepackage{abbreviazioni}

\begin{document}

\title{On Bohr's equivalence theorem}

\author{Mattia Righetti}

\address{Dipartimento di Matematica, Universit\`a di Genova, via Dodecaneso 35, 16146, Genova, Italy}
\curraddr{Centre de Recherches Math\'ematiques, Universit\'e de Montr\'eal, P.O. box 6128, Centre-Ville Station, Montr\'eal, Qu\'ebec H3C 3J7, Canada.}
\email{righetti@dima.unige.it}           

\date{}

\begin{abstract}
In this note we prove a converse of Bohr's equivalence theorem for Dirichlet series under some natural assumptions.\\

\end{abstract}

\maketitle

\renewcommand*{\thetheorem}{\Alph{theorem}}

\input{intro}
\begin{acknowledgements}This paper is part of my Ph.D. thesis at the Department of Mathematics of the University of Genova. I express my sincere gratitude to my supervisor Professor Alberto Perelli for his support and for the many valuable suggestions. I wish to thank also Michel Balazard for the above remark on the sequence of exponents.
\end{acknowledgements}
\input{proof_converse}

\bibliographystyle{amsplain}      
\bibliography{biblio}   

\end{document}

%% file: intro.tex
\section{Introduction}

Bohr's interest in the Riemann zeta function led him to study the set of values taken by Dirichlet series in their half-plane of absolute convergence. For this problem Bohr developed a new method: associating to any Dirichlet series a power series with infinitely many variables (see \cite{bohr10}). He then introduced an equivalence relation among Dirichlet series and showed that equivalent Dirichlet series take the same set of values in certain open half-planes (see \cite{bohr7}). We give here a very brief account of this theory; for a complete treatment we refer to Bohr's original work \cite{bohr7} and to Chapter 8 of Apostol \cite{apostol}.

We call \emph{general Dirichlet series} any complex function $f(s)$, in the variable $s=\sigma+it$, that has a series representation of the form
$$f(s)=\sum_{n=1}^\infty a(n)\exp{-\lambda(n)s},$$
where the \emph{coefficients} $a(n)$ are complex and the sequence of \emph{exponents} $\Lambda=\{\lambda(n)\}$ consists of real numbers such that $\lambda(1)<\lambda(2)<\cdots$ and $\lambda(n)\rightarrow\infty$ as $n\rightarrow\infty$.\\
Note that this class of general Dirichlet series includes both power series, when $\lambda(n)=n$, and ordinary Dirichlet series, when $\lambda(n)=\log(n)$.

\begin{remark*} The above definition of general Dirichlet series is the one that is given in the work of Bohr \cite{bohr7}, and it is more restrictive then the usual definition without conditions on the sequence of exponents, which is already present in later works of Bohr (see e.g. \cite{bohr5}). Restricting to the above setting has the advantage that the region of absolute convergence of the series is a right half-plane, like for ordinary Dirichlet series (see e.g. \cite[\S8.2]{apostol}). Since this is fundamental in the proof of the converse theorem, we have decided to work in the same setting of \cite{bohr7}. However, as we remarked in our Ph.D. thesis \cite{righettiphd}, Bohr's equivalence theorem holds true \emph{mutatis mutandis} even for the usual general Dirichlet series.
\end{remark*}

Following Bohr (see e.g. \cite[\S8.3]{apostol}), given a sequence of exponents $\Lambda=\{\lambda(n)\}$, we say that a sequence of real numbers $B=\{\beta(n)\}$ is a \emph{basis} for $\Lambda$ if it satisfies the following conditions:
\begin{enumerate}[(i)]
\item the elements of $B$ are linearly independent over the rationals;
\item for every $n$, $\lambda(n)$ is expressible as a finite linear combination over $\Q$ of elements of $B$;
\item for every $n$, $\beta(n)$ is expressible as a finite linear combination over $\Q$ of elements of $\Lambda$.
\end{enumerate}
We may express the above conditions in matrix notation by considering $\Lambda$ and $B$ as infinite column vectors (see \cite[\S8.4]{apostol}). In particular, if $B$ is a basis for $\Lambda$, we may write $\Lambda=RB$ and $B=T\Lambda$ for some \emph{Bohr matrices} $R$ and $T$.

We fix a sequence of exponents $\Lambda =\{\lambda(n)\}$ and a basis $B$ of $\Lambda$, so that we may write $\Lambda=R B$. Consider two general Dirichlet series with the same sequence of exponents $\Lambda$, say
$$f(s)=\sum_{n=1}^\infty a(n)\exp{-\lambda(n)s}\quad\hbox{and}\quad g(s)=\sum_{n=1}^\infty b(n)\exp{-\lambda(n)s}.$$
Then, we say that $f(s)$ and $g(s)$ are \emph{equivalent} ($f\sim g$), with respect to $B$, if there exists a sequence of real numbers $Y=\{y(n)\}$ such that
\begin{equation*}\label{eq:equiv}
b(n)=a(n)\exp{i(RY)_n}\quad\hbox{ for every }n.
\end{equation*}

We may now state Bohr's equivalence theorem  (cf. Theorem 8.16 of Apostol \cite{apostol}), which is, roughly speaking, a combination of Kronecker's approximation theorem, Rouch\'e's theorem and the absolute convergence of the Dirichlet series. 

\begin{theorem}[Bohr, {\cite[Satz 4]{bohr7}}]\label{theorem:bohr_equivalence}
Let $f_1(s)$ and $f_2(s)$ be equivalent general Dirichlet series absolutely convergent for $\sigma>\alpha$. Then in any open half-plane $\sigma>\sigma_0\geq \alpha$ the functions $f_1(s)$ and $f_2(s)$ take the same set of values.
\end{theorem}

\begin{remark*} Although Bohr's equivalence theorem is usually stated for open half-planes, we have already remarked in \cite{righetti} and \cite{righettiphd} that Theorem \ref{theorem:bohr_equivalence} holds true also for open vertical strips.
\end{remark*}

In particular, one gets immediately the following more practical version of the above theorem, in the sense that in the applications to the value distribution of $L$-functions one usually appeals to this statement. Similar results in particular cases may be found for example in Bohr \cite[\S2]{bohr8}, Titchmarsh \cite[\S11.4]{titchmarsh}, Bombieri and Mueller \cite[Lemma 1]{bombierimueller} and Bombieri and Ghosh \cite[p. 240]{bombierighosh}

\begin{theorem}[Bohr]\label{theorem:practical_bohr}
Let $f(s)$ be a general Dirichlet series absolutely convergent for $\sigma>\alpha$, and let $S_f(\sigma_1,\sigma_2)$ is the set of values taken by $f(s)$ in the strip $\alpha\leq\sigma_1<\sigma<\sigma_2\leq \infty$, then
$$S_f(\sigma_1,\sigma_2) = \{g(\sigma)\mid \sigma\in(\sigma_1,\sigma_2),\,g \sim f \}.$$
If moreover $\Lambda$ has an \emph{integral basis} B, i.e. the Bohr matrix $R$ such that $\Lambda=RB$ has only integer entries, and $V_f(\sigma_0)$ is the set of values taken by $f(s)$ on the vertical line $\sigma=\sigma_0>\alpha$, then
\begin{equation}\label{eq:equiv_local}
\conj{V_f(\sigma_0)}=\{g(\sigma_0)\mid g\sim f\}.
\end{equation}
\end{theorem}
The first statement follows immediately from Theorem \ref{theorem:bohr_equivalence}, while the second is Satz 3 in \cite{bohr7}.

\begin{remark*} Note that \eqref{eq:equiv_local} doesn't hold in general. Indeed, consider the following example given by Bohr \cite[pp. 151--153]{bohr7}: $\lambda(n)=2n-1+\frac{1}{2(2n-1)}$, $f(s) = \sum_n \exp{-\lambda(n)s}$ and $g(s)=-f(s)$. Then the only equivalent functions to $f(s)$ are its vertical shifts and, as proved by Bohr, $g(s)$ is not equivalent to $f(s)$. On the other hand if $\tau_m = 2\pi \prod_{n\leq m} (2n-1)$ then $f(s+i\tau_m)$ converges uniformly to $g(s)$, and conversely $g(s+i\tau_m)$ converges uniformly to $f(s)$. Hence, using Rouch\'e's theorem one may show that $g(\sigma_0)\in \conj{V_f(\sigma_0)}$ for any $\sigma_0>\alpha$.

Actually the same argument shows that $f(s)$ and $g(s)$ have the same values in any open right half-plane, so in general there is no converse to Theorem \ref{theorem:bohr_equivalence}. This is because, as the example shows, the set of general Dirichlet series equivalent to a certain general Dirichlet series may not be closed with respect to uniform convergence on compact subsets. However this is a closed set if the sequence of exponents has an integral basis, e.g. for ordinary Dirichlet series.
\end{remark*}

This leads us to the following converse theorem.

\begin{theorem}\label{theorem:converse_bohr} Let $f_1(s)$ and $f_2(s)$ be general Dirichlet series with the same sequence of exponents $\Lambda=\{\lambda(n)\}$ and absolutely converging for $\sigma>\alpha$ for some $\alpha$. Suppose that $f_1(s)$ and $f_2(s)$ take the same set of values in any open half-plane $\sigma>\sigma_0> \alpha$. Then $f_2(s)$ belongs to the closure, with respect to uniform convergence on compact susbsets of $\sigma>\alpha$, of the set of general Dirichlet series equivalent to $f_1(s)$ and vice versa.\\
If furthermore $\Lambda$ has an integral basis, then $f_1(s)$ and $f_2(s)$ are equivalent.
\end{theorem}

\begin{remark*} The requirement that $f_1(s)$ and $f_2(s)$ should have the same exponents is not really restrictive: one just needs to take the union of the sets of exponents. On the other hand \textit{a fortiori} they would have the same exponents in the sense that $a_1(n)=0$ if and only if $a_2(n)=0$.
\end{remark*}

%% file: proof_converse.tex
\section{Proof of Theorem \ref{theorem:converse_bohr}}

Let $B=\{\beta(j)\}$ be a basis for $\Lambda$ and let $R=(r_{ij})$ be the Bohr matrix such that $\Lambda=RB$. We work by induction: we want to show that for any $n$ there exist $Y_n$ such that
\begin{equation}\label{eq:induction}
a_2(k)=a_1(k)\exp{i(RY_n)_k}, \qquad k=1,\ldots, n. 
\end{equation}

For $n=1$, we show that there exist sequences $\{\sigma_m\}$, with $\sigma_m\rightarrow\infty$, $\{t_m\}$ and $\{t_m'\}$ such that
\begin{equation}\label{eq:limit_1}
\exp{\lambda(1)\sigma_m}[f_1(\sigma_m+it_m)-f_2(\sigma_m+it_m')]\rightarrow 0 , \qquad m\rightarrow\infty.
\end{equation}
We construct $\{\sigma_m\}$ in the following way. It is well known that for every $v\in\C$ the Dirichlet series has a zero-free right half-plane which is maximal, say $\sigma>\sigma^*(v)$. Hence we take $\sigma_m = \sigma^*(f_2(m))$, $m=1,2,\ldots$. Note that by definition and by hypothesis for any $\eps>0$ we have
$$f_2(m)\in S_{f_j}(\sigma_m-\eps,\infty)\backslash S_{f_j}(\sigma_m,\infty),\qquad j=1,2.$$
Therefore, for any $m$, there exist $\sigma_{m,1}$ and $\sigma_{m,2}$ such that $0<\sigma_m-\sigma_{m,j}<\eps$, $j=1,2$, and $t_m$ and $t'_m$ such that
$$f_1(\sigma_{m,1}+it_m)=f_2(\sigma_{m,2}+it_m')=f_2(m).$$
Hence, taking $\eps=\exp{-|\lambda(1)|\sigma_m}/(|\lambda(1)|+1)$ we get
\begin{spliteq*}
|f_1(\sigma_m+it_m)-f_2(\sigma_m+it_m')| &\leq |f_1(\sigma_m+it_m)-f_1(\sigma_{m,1}+it_m)|\\
&\qquad\qquad+|f_2(\sigma_{m,2}+it_m')-f_2(\sigma_m+it_m')| \\
& \leq \sum_{n=1}^\infty |a_1(n)|\exp{-\lambda(n)\sigma_{m}}| \exp{\lambda(n)(\sigma_m-\sigma_{m,1})}- 1|\\
&\qquad +\sum_{n=1}^\infty |a_2(n)|  \exp{-\lambda(n)\sigma_{m}}|\exp{\lambda(n)(\sigma_m-\sigma_{m,2})}-1|\\
& \ll (|a_1(1)|+|a_2(1)|)\exp{-(\lambda(1)+|\lambda(1)|)\sigma_m}\\
&\qquad\qquad+ \sum_{n=2}^\infty |a_1(n)|\exp{-\lambda(n)\sigma_{m,1}}+\sum_{n=2}^\infty|a_2(n)|\exp{-\lambda(n)\sigma_{m,2}}.
\end{spliteq*}
From this we have that \eqref{eq:limit_1} immediately follow.

On the other hand, we recall that $R$ has rational entries, so we may write $r_{ij}=b_{ij}/q_{ij}$, and for any $i$ we have $r_{ij}\neq 0$ only for a finite number $j$s. Hence we define 
$$d_h = \hbox{l.c.m.} \{q_{ij}\: :\: i=1,\ldots,h,\, j\hbox{ s.t. } r_{ij}\neq0\},\qquad h=1,2,\ldots.$$
If we apply Helly's selection principle (see e.g. Lemma 1, \S8.12 of \cite{apostol}) to the uniformly bounded double sequences $\{-t_m\beta(j)$ mod $2\pi d_1\}$ and $\{-t_m'\beta(j)$ mod $2\pi d_1\}$, then there exist sequences $\Theta=\{\theta(j)\}$ and $\Theta'=\{\theta'(j)\}$ and a subsequence $\{m_k\}$ such that, when $k\rightarrow\infty$,
$$\exp{\lambda(1)\sigma_{m_k}}[f_1(\sigma_{m_k}+it_{m_k})-f_2(\sigma_{m_k}+it_{m_k}')]\rightarrow a_1(1)\exp{i (R\Theta)_1}-a_2(1)\exp{i (R\Theta')_1}.$$
By \eqref{eq:limit_1} we have that \eqref{eq:induction} holds for $n=1$ by taking $Y_1=\Theta-\Theta'$.

Suppose now that \eqref{eq:induction} holds, then reasoning as above it is easy to get that \eqref{eq:induction} holds also for $n+1$. Indeed, let $\{\sigma_m\}$ be the same sequence as above and let $g_n(s)$ be the general Dirichlet series equivalent to $f_1(s)$ defined by
$$g_n(s) = \sum_k a_1(k)\exp{i(RY_n)_k}\exp{-\lambda(k)s}.$$
By hypothesis and Theorem \ref{theorem:bohr_equivalence} we have that for any $\eps>0$ $$f_2(m)\in S_{g_n}(\sigma_m-\eps,\infty)\backslash S_{g_n}(\sigma_m,\infty).$$
Therefore, for any $m$, there exist $\sigma_{m,1}$ and $\sigma_{m,2}$ such that $0<\sigma_m-\sigma_{m,j}<\eps$, $j=1,2$, and $t_m$ and $t'_m$ such that
$$g_n(\sigma_{m,1}+it_m)=f_2(\sigma_{m,2}+it_m')=f_2(m).$$
Hence, taking $\eps=\exp{-|\lambda(n+1)|\sigma_m}/(|\lambda(n+1)|+1)$, since \eqref{eq:induction} holds, analogously as before we get
\begin{spliteq*}
|g_n(\sigma_m+it_m)-f_2(\sigma_m+it_m')| 
& \ll  (|a_1(n+1)|+|a_2(n+1)|)\exp{-(\lambda(n+1)+|\lambda(n+1)|)\sigma_m}\\
&\qquad\qquad+ \sum_{k=n+2}^\infty |a_1(k)|\exp{-\lambda(k)\sigma_{m,1}}+\sum_{k=n+2}^\infty|a_2(k)|\exp{-\lambda(k)\sigma_{m,2}}
\end{spliteq*}
From this we deduce that 
\begin{equation*}
\exp{\lambda(n+1)\sigma_m}[g_n(\sigma_m+it_m)-f_2(\sigma_m+it_m')]\rightarrow 0 , \qquad m\rightarrow\infty.
\end{equation*}

On the other hand, as before, if we apply Helly's selection principle to the uniformly bounded double sequences $\{-t_m\beta(j)$ mod $2\pi d_{n+1}\}$ and $\{-t_m'\beta(j)$ mod $2\pi d_{n+1}\}$, then there exist $\Theta=\{\theta(j)\}$ and $\Theta'=\{\theta'(j)\}$ and a subsequence $\{m_k\}$ such that, when $k\rightarrow\infty$,
$$\exp{\lambda(n+1)\sigma_{m_k}}[g_n(\sigma_{m_k}+it_{m_k})-f_2(\sigma_{m_k}+it_{m_k}')]\rightarrow a_1(n+1)\exp{i (R(Y_n+\Theta))_{n+1}}-a_2(n+1)\exp{i (R\Theta')_{n+1}}.$$
By the uniqueness of the limit we have that \eqref{eq:induction} holds for $n+1$ by taking $Y_{n+1}=Y_n+\Theta-\Theta'$.

Finally we note that by induction we have actually constructed a sequence $\{g_n(s)\}$ of general Dirichlet series equivalent to $f_1(s)$ and by \eqref{eq:induction} we have that $g_n(s)$ converges uniformly on every compact subset of $\sigma>\alpha$ to $f_2(s)$. Since we may change $f_1(s)$ with $f_2(s)$, this proves the first statement.

If $B$ is an integral basis for $\Lambda$, then we may apply Helly's selection principle to the uniformly bounded double sequence $\{y_n(j)$ mod $2\pi\}$, where $Y_n=\{y_n(j)\}$ is the sequence obtained at each step of the induction process. Hence there exist a sequence $Y$ and a subsequence $\{n_h\}$ such that for $h\rightarrow\infty$ the functions $g_{n_h}(s)$ converge uniformly on every compact subset of $\sigma>\alpha$ to 
\begin{equation*}
g(s) = \sum_k a_1(k)\exp{i(RY)_k}\exp{-\lambda(k)s},
\end{equation*}
which is equivalent to $f_1(s)$. By the uniqueness of the limit and of the analytic continuation we have that $g(s)=f_2(s)$ in the half-plane of absolute convergence, and the result follows.